\documentclass[10pt, titlepage]{amsart}

\usepackage{ae} 
\usepackage[T1]{fontenc}
\usepackage[cp1250]{inputenc}
\usepackage{amsmath}
\usepackage{amssymb, amsfonts,amscd,verbatim}
\usepackage[dvips]{graphicx}
\usepackage{latexsym}
\usepackage{indentfirst}
\usepackage{latexsym}

\usepackage{graphicx}
\usepackage{verbatim}   
\usepackage{color}      
\usepackage{subfigure}  

\usepackage{amsmath}    

\theoremstyle{plain}
\newtheorem{Prop}{Proposition}[section]
\newtheorem{Thm}[Prop]{Theorem}
\newtheorem{Cor}[Prop]{Corollary}

\theoremstyle{definition}
\newtheorem{Def}[Prop]{Definition}

\theoremstyle{remark}

\newtheorem{Example}[Prop]{\bf Example}

\def\int{\mathop{\roman{int}}}

\def\1{^{-1}}

\def\p{{\mathbf p}}
\def\q{{\mathbf q}}

\def\asdim{\mathrm{asdim}}

\def\dokaz{{\bf Proof. }}
\def\edokaz{\hfill $\blacksquare$}

\def\RipsG{\mathrm{RipsG}}
\def\Rips{\mathrm{Rips}}

\numberwithin{equation}{section}

\input pstricks.tex
\input xy
\xyoption{all}

\begin{document}
\title[
A combinatorial approach to coarse geometry
]%
   {A combinatorial approach to coarse geometry}

\author{M.~Cencelj}
\address{IMFM,
Univerza v Ljubljani,
Jadranska ulica 19,
SI-1111 Ljubljana,
Slovenija }
\email{matija.cencelj@guest.arnes.si}

\author{J.~Dydak}
\address{University of Tennessee, Knoxville, TN 37996, USA}
\email{dydak@math.utk.edu}

\author{A.~Vavpeti\v c}
\address{Fakulteta za Matematiko in Fiziko,
Univerza v Ljubljani,
Jadranska ulica 19,
SI-1111 Ljubljana,
Slovenija }
\email{ales.vavpetic@fmf.uni-lj.si}

\author{\v Z.~Virk}
\address{University of Tennessee, Knoxville, TN 37996, USA}
\email{zigavirk@gmail.com}

\date{ \today
}
\keywords{asymptotic dimension, coarse geometry, simplicial trees}

\subjclass[2000]{Primary 54F45; Secondary 55M10}

\thanks{This research was supported by the Slovenian Research
Agency grants P1-0292-0101 and J1-2057-0101.}
\thanks{The second-named author was partially supported
by MEC, MTM2006-0825.}

\begin{abstract}

Using ideas from shape theory we embed the coarse category of metric spaces
into the category of direct sequences of simplicial complexes with bonding maps
being simplicial. Two direct sequences of simplicial complexes are equivalent if one of them can be transformed
to the other by contiguous factorizations of bonding maps and by taking infinite subsequences.
That embedding can be realized by either Rips complexes or analogs of Roe's
anti-\v Cech approximations of spaces.

In that model coarse $n$-connectedness of $\mathcal{K}=\{K_1\to K_2\to\ldots\}$ means that for each $k$ there is $m > k$
such that the bonding map from $K_k$ to $K_m$ induces trivial homomorphisms of all homotopy groups up to and including $n$.

The asymptotic dimension being at most $n$ means that for each $k$ there is $m > k$
such that the bonding map from $K_k$ to $K_m$ factors (up to contiguity)
through an $n$-dimensional complex.

Property A of G.Yu is equivalent to the condition that for each $k$ and for each $\epsilon > 0$ there is $m > k$
such that the bonding map from $\vert K_k\vert $ to $\vert K_m\vert$ 
has a contiguous approximation $g\colon \vert K_k\vert \to\vert K_m\vert$
which sends simplices of $\vert K_k\vert$ to sets of diameter
at most $\epsilon$. \end{abstract}

\maketitle

\medskip
\medskip
\tableofcontents
\section{Introduction}

In homotopy theory the class with optimal properties consists of CW complexes.
Other spaces are investigated by mapping CW complexes to them. That leads
to the concept of the {\bf singular complex} $Sin(X)$ of a space $X$
together with the projection $p\colon Sin(X)\to X$ so that the following
conditions are satisfied:
\begin{enumerate}
\item Any continuous map $f\colon K\to X$, $K$ a CW complex, lifts up to homotopy
to a continuous map $g\colon K\to Sin(X)$ and $g$ is unique up to homotopy,
\item $p$ induces isomorphisms of all singular homology groups.
\end{enumerate}

Thus $p\colon Sin(X)\to X$ is a universal object among all continuous maps from CW complexes to $X$ and
$Sin(X)$ represents all information about $X$ from the point of view of the weak homotopy theory. A map $f\colon X\to Y$ of topological spaces is a {\bf weak homotopy equivalence}
if it induces a bijection $f_\ast\colon [K,X]\to [K,Y]$ of sets of homotopy classes
for all CW complexes $K$ (equivalently, it induces a homotopy equivalence from $Sin(X)$ to $Sin(Y)$).

The shape theory (see \cite{DydSeg}
and \cite{MarSeg} ) represents the dual point of view:
an arbitrary topological space $X$ is investigated by mapping $X$ to CW complexes $K$.
A map $f\colon X\to Y$ of topological spaces is a {\bf shape equivalence}
if it induces a bijection $f^\ast\colon [Y,K]\to [X,K]$ of sets of homotopy classes
for all CW complexes $K$.

In contrast to the weak homotopy theory there is no universal
continuous map from $X$ to a particular CW complex. Instead, each space $X$ has
the {\bf \v Cech system} $\{K_\alpha,[p^\beta_\alpha],A\}$
with projections $p_\alpha\colon X\to K_\alpha$
which reflects the shape of $X$ in the following sense:
\begin{itemize}
\item [a.] For any continuous map $f\colon X\to K$ to a CW complex $K$ there is $\alpha\in A$
and $g\colon K_\alpha\to K$ such that $f$ is homotopic to $g\circ p_\alpha$,
\item[b.] Given $\alpha\in A$ and given continuous maps $g,h\colon K_\alpha\to K$
with $g\circ p_\alpha $ homotopic to $h\circ p_\alpha$
there is $\beta \ge \alpha$ so that $g\circ p^\beta_\alpha$ is homotopic
to $h\circ p^\beta_\alpha$.
\end{itemize}

The \v Cech system of $X$ consists of geometric realizations of nerves of numerable open coverings 
of $X$ and the bonding maps are simplicial maps induced by functions between covers
that reflect one of them being a refinement of the other.

In this paper we will show that the coarse category
of metric spaces is dual to shape category in the following sense:
for each $X$ we consider the set of equivalence
classes $C(K,X)$ of functions $f\colon K\to X$ from
a simplicial complex $K$ to $X$.
It is required that the family $\{f(\Delta)\}_{\Delta\in K}$ is uniformly bounded
(we call such functions {\bf bornologous}) and $f$ and $g$ are equivalent if the family $\{f(\Delta)\cup g(\Delta)\}_{\Delta\in K}$ is uniformly bounded.
Each metric space has a coarse \v Cech system
$\{K_m\}$ consisting of a direct system of simplicial complexes and simplicial maps.
That system is universal among all bornologous functions from
simplicial complexes to $X$.

We show asymptotic dimension is dual to the shape
dimension and coarse connectivity is dual to shape connectivity.

{\bf Convention}: The metrics considered in this paper may attain infinite values.
As explained in \cite{DH}, the primary advantage of such metrics is the ability of constructing disjoint unions of metric spaces (see the proof of \ref{CharOfBorno}).

We will make a careful distinction between graphs (or simplicial complexes) and
their geometric realization.
Usually we will focus on the set of vertices $X=G^{(0)}$ of a graph $G$
(or a simplicial complex $K$) and we may denote $G$ (or $K$) by $X$ in the absence
of other graph (or simplicial complex) structures on $X$.
By {\bf geometric realization} $|G|$ (or $|K|$) we mean $G^{(0)}$ and the union
of all geometric edges (geometric simplices) induced by edges in $G$ (simplices in $K$).
Thus, given a simplex $\Delta=[x_0,\ldots,x_n]$ (a finite subset of $X$),
its geometric realization $|\Delta|$ is the set of all formal
linear combinations $\sum\limits_{i=0}^n t_i\cdot x_i$,
where $t_i\ge 0$ and $\sum\limits_{i=0}^n t_i=1$.

Given a set $X$ we can put a graph structure $G$ on it by specifying
all the edges (example: the Cayley graph of a finitely generated group).
A graph structure $G$ on $X$ leads to the {\bf graph metric} $d_G$ on $X$ as follows
(notice the advantage of using metrics with infinite values):
 \begin{enumerate}
\item $d_G(v,w)=0$ if and only if $v=w$,
\item $0 < d_G(v,w)\leq n$ if and only if $v\ne w$
and there is a chain $v_0=v,\ldots,v_n=w$
such that $[v_i,v_{i+1}]$ is an edge in $G$ for each
$0\leq i < n$. $d_G(v,w)=k > 0$ if $d_G(v,w)\leq k$ and $d_G(v,w)\leq k-1$ is false,
\item $d_G(v,w)=\infty$ otherwise.
\end{enumerate}

Conversely, given a metric $d_X$ on $X$, we can put several graph structures
on $X$:
\begin{Def}
Given a metric space $(X,d_X)$ and $t > 0$,
{\bf the Rips graph} $\RipsG_t(X)$ consists of edges
$[x,y]$ such that $d_X(x,y) \leq t$.
\end{Def}

Observe $\RipsG_n(\RipsG_m(X,d_X))=\RipsG_{m\cdot n}(X,d_X)$
for $(X,d_X)$ geodesic and $m,n$ positive integers.

Notice that $\pi_t\colon \RipsG_t(X)\to X$ induced
by the identity function is $t$-Lipschitz.
Moreover, $t$-Lipschitz maps from a graph $V$
to $X$ are in one-to-one correspondence
with short maps from $V$ to $\RipsG_t(X)$
(following Gromov
by {\bf short maps} we mean Lipschitz maps with the Lipschitz constant one).

Notice that a function $f \colon G\to H$ between graphs is short if and only if
for every edge $[x,y]$ of G,  $[f(x),f(y)]$ is an edge of $H$.

Observe that the metric $d_X$ of a metric space $(X,d_X)$ is induced by
some graph structure on $X$ if and only if $d_X$  is integer valued and $X$ is
$1$-geodesic. A metric space $(X,d_X)$ is {\bf $t$-geodesic}
if for every two points $x$ and $y$ there is a $t$-chain $x_0=x$, \ldots, $x_k=y$
(that means $d_X(x_i,x_{i+1})\leq t$ for all $0\leq i\leq k-1$)
such that $d_X(x,y)=\sum\limits_{i=0}^{k-1}d_X(x_i,x_{i+1})$.
Furthermore the graph metric $d_G$ on $G$ can be extended to
a geodesic metric on the geometric realization $|G|$ of $G$.

One can generalize the concept of Rips graphs from metric spaces to
arbitrary sets provided a cover of the set is given:

\begin{Def}
Given a set $X$ and a cover $\mathcal{U}$ of $X$,
{\bf the Rips graph} $\RipsG_{\mathcal{U}}(X)$ consists of edges
$[x,y]$ such that there is an element $U$ of $\mathcal{U}$ containing $x$ and $y$.
\end{Def}

Notice $\RipsG_t(X)\subset \RipsG_{\mathcal{U}}(X)\subset \RipsG_{2t}(X)$,
where $\mathcal{U}$ is the cover of $X$ by closed $t$-balls.

\section{The coarse category}

In this section we introduce the coarse category of metric spaces 
in a way different (but equivalent) to that of Roe \cite{Roe lectures}. It is similar
to the intuitive way of introducing asymptotic cones of metric spaces: one
looks at a given metric space from farther and farther away. In our case the basic
idea is to look at two functions from farther and farther away.

From the large scale point of view two functions $f,g\colon (X,d_X)\to (Y,d_Y)$
should be considered indistinguishable if they are within finite distance from
each other (i.e., $\sup\limits_{x\in X} d_Y(f(x),g(x)) < \infty$).
Therefore it makes sense to consider the set of equivalence classes $[X,Y]_{ls}$
induced on the set of all functions from $X$ to $Y$ by the equivalence relation $f\sim_{ls} g$ if and only if
$\sup\limits_{x\in X} d_Y(f(x),g(x)) < \infty$.
\par
The next logical step, from the categorical point of view, is to consider functions that preserve the relation $\sim_{ls}$. Obviously,
$f\sim_{ls}g$ implies $f\circ\alpha\sim_{ls}g\circ\alpha$ for any $\alpha\colon  (Z,d_Z)\to (X,d_X)$, so the following definition addresses the crux of the matter:

\begin{Def}
A function $\alpha\colon (X,d_X)\to (Y,d_Y)$ of metric spaces is {\bf bornologous}
(or {\bf large scale uniform}) if $f\sim_{ls}g$, $f,g\colon (Z,d_Z)\to (X,d_X)$,
implies $\alpha \circ f \sim_{ls}\alpha\circ g$.
\end{Def}

Let us show that our definition of a function being bornologous is equivalent to that of J.Roe \cite{Roe lectures}:
\begin{Prop}\label{CharOfBorno}
If $\alpha\colon (X,d_X)\to (Y,d_Y)$ is a function of metric spaces, then the following
conditions are equivalent:
\begin{itemize}
\item[a.] $\alpha$ is bornologous,
\item[b.] For any $t > 0$ there is $s > 0$ such that
$d_X(x,y) < t$ implies $d_Y(\alpha(x),\alpha(y)) < s$.
\end{itemize}

\end{Prop}
\dokaz Suppose $\alpha$ is bornologous. For $t>0$ define a metric space
$X_t=\coprod_{x\in X} B(x,t)\times \{x\}$ with metric
$$
d_t((x_1,y_1),(x_2,y_2))=\left\{\begin{array}{ll}
d_X(x_1,x_2) &; y_1=y_2,\\
\infty &; y_1\ne y_2
\end{array}\right. .
$$
Maps $f,g\colon X_t\to X$ defined as $f(x,y)=x$ and $g(x,y)=y$ are at distance at most $t$, hence $f \sim_{ls} g$.
Because $\alpha$ is bornologous, $\alpha \circ f \sim_{ls}\alpha\circ g$ which means that there exists
$s>0$ such that $d_X((\alpha\circ f)(x,y),(\alpha\circ g)(x,y))<s$. Let $d_X(x,y)<t$,
then $(x,y)\in X_t$. Because $(\alpha\circ f)(x,y)=\alpha(x)$ and $(\alpha\circ g)(x,y)=\alpha(y)$,
the distance $d_Y(\alpha(x),\alpha(y))<s$.

Suppose that for any $t > 0$ there is $s > 0$ such that
$d_X(x,y) < t$ implies $d_Y(\alpha(x),\alpha(y)) < s$.
Let $f,g\colon (Z,d_Z)\to (X,d_X)$ satisfy $\sup_{z\in Z} d_X(f(x),g(x))=t<\infty$. Let $s$ be as above.
Then $d_X(\alpha(f(z)),\alpha(g(z)))<s$, hence $\alpha \circ f \sim_{ls}\alpha\circ g$.

\edokaz

\begin{Def}
A bornologous function $\alpha\colon (X,d_X)\to (Y,d_Y)$ of metric spaces is a {\bf large scale isomorphism} if it induces a bijection $\alpha_\ast\colon [Z,X]_{ls}\to [Z,Y]_{ls}$
for all metric spaces $(Z,d_Z)$.
\end{Def}

As in the case of bornologous functions, our definition
of large scale isomorphisms is equivalent to that in \cite{Roe lectures}:

\begin{Prop}
If $\alpha\colon (X,d_X)\to (Y,d_Y)$ is a bornologous function of metric spaces, then the following
conditions are equivalent:
\begin{itemize}
\item[a.] $\alpha$ is a large scale isomorphism,
\item[b.] There is a bornologous function $\beta\colon (Y,d_Y)\to (X,d_X)$
such that $\alpha\circ\beta\sim_{ls} id_Y$ and $\beta\circ\alpha\sim_{ls} id_X$.
\end{itemize}

\end{Prop}
\dokaz Suppose $\alpha$ is a large scale isomorphism. The map $\alpha_*\colon [Y,X]_{ls}\to [X,X]_{ls}$ is a bijection, hence
there exists a map $\beta\colon Y\to X$ such that $\alpha\circ\beta\sim_{ls} id_Y$. Because $\alpha\circ\beta\circ\alpha\sim_{ls}\alpha\circ id_X$
and $\alpha_*\colon [X,X]_{ls}\to [X,Y]_{ls}$ is a bijection, also $\beta\circ\alpha\sim_{ls} id_X$.
Let us show that $\beta$ is bornologous. Let $f,g\colon Z\to Y$ be maps between metric spaces and $f\sim_{ls} g$. Because $\alpha\circ \beta$
is bornologous, $(\alpha\circ\beta)\circ f\sim_{ls} (\alpha\circ \beta)\circ g$. Because $\alpha_*\colon [Z,X]_{ls} \to [Z,Y]_{ls}$
is a bijection, $\beta\circ f\sim_{ls} \beta\circ g$.

Suppose there is a bornologous function $\beta\colon (Y,d_Y)\to (X,d_X)$
such that $\alpha\circ\beta\sim_{ls} id_Y$ and $\beta\circ\alpha\sim_{ls} id_X$.
Let $Z$ be a metric space and $f,g\colon Z\to X$ such that $\alpha\circ f\sim_{ls} \alpha\circ g$.
Then $f\sim_{ls} \beta\circ\alpha\circ f\sim_{ls} \beta\circ\alpha\circ g\sim_{ls} g$, hence $\alpha_*\colon [Z,X]_{ls}\to [Z,Y]_{ls}$
is a monomorphism. Let $h\colon Z\to Y$ be a map. Then $\alpha\circ (\beta\circ h)\sim_{ls} h$, so $\alpha_*$ is an epimorphism.
\edokaz

In view of \ref{CharOfBorno} the simplest bornologous functions
are Lipschitz functions. Conversely, it is easy to deduce from
 \ref{CharOfBorno} that every bornologous function
 defined on a $t$-geodesic space is Lipschitz.

The following result shows that Lipschitz maps from graphs are of primary interest
in large scale geometry.

\begin{Thm}
Let $f\colon (X,d_X)\to (Y,d_Y)$ 
be a function of metric spaces.
\begin{itemize}
\item [a.] $f$ is bornologous if and only if for every Lipschitz
function $g\colon (G,d_G)\to X$, $G$ any graph, $f\circ g$ is Lipschitz,
\item[b.] $f$ induces a large scale isomorphism if and only
if every Lipschitz function $h\colon G\to Y$, $G$ any graph,
lifts (up to large scale equivalence) to a Lipschitz map to $X$ and the lift is unique up to large scale equivalence.
\end{itemize}

\end{Thm}
\dokaz (a).
Let a map $f\colon (X,d_X)\to (Y,d_Y)$ be bornologous and $g\colon (G,d_G)\to X$ be a Lipschitz map with
Lipschitz constant $t$. There exists $s>0$ such that $d_X(x,y)<t$ implies $d_Y(f(x),f(y))<s$.
Let $a,b\in G$ and $d_G(a,b)=n$, then there exit a sequence $a_0,\ldots, a_n\in G$ such that $a_0=a$, $a_n=b$,
and $d_G(a_{i-1},a_i)=1$ for $i=1,\ldots, n$. For $i=1,\ldots, n$ the distance $d_X(g(a_{i-1}),g(a_i))<t$, hence $d_Y(fg(a_{i-1}),fg(a_i))<s$.
Therefore
$$
d_Y(fg(a),fg(b))\le \sum_{i=1}^n d_Y(fg(a_{i-1}),fg(a_i))\le \sum_{i=1}^n s=s d_G(a,b)
$$
and the map $f\circ g$ is Lipschitz.

Suppose that for every Lipschitz function $g\colon (G,d_G)\to X$, $G$ any graph,
$f\circ g$ is Lipschitz.
Let us prove that $f$ is bornologous. Let $t>0$. The map $\pi_t\colon G= \RipsG_t(X)\to X$ is $t$-Lipschitz,
hence the map $f\circ\pi_t$ is $L$-Lipschitz for some $L$. If $d_X(x,y) < t$ and $x\ne y$, then
$$
d_Y(f(x),f(y))=d_Y(f\pi_t(x),f\pi_t(y))\le Ld_G(x,y)=L,
$$
so $f$ is bornologous.

(b).
Let $f$ be a large scale isomorphism. There exists a bornologous function $g\colon Y\to X$, such that $f\circ g\sim_{ls} id_Y$
and $g\circ f\sim_{ls} id_X$. Let $h\colon G\to Y$ be a Lipschitz map. Let $\widetilde h=g\circ h$. Then $f\circ \widetilde h=f\circ g\circ h
\sim_{ls} h$, so $\widetilde h$ is a lift up to large scale equivalence of the map $h$. Suppose $h'\colon G\to X$ be a Lipschitz map
such that $f\circ h'\sim_{ls} h$. Then $\widetilde h=g\circ h\sim_{ls} g\circ f\circ h'\sim_{ls} h'$, hence
a lift is unique.

\par Suppose every Lipschitz function $h\colon G\to Y$, $G$ any graph,
lifts (up to large scale equivalence) to $X$ and the lift is unique up to large scale equivalence.
Consider lifts $h_t\colon \RipsG_t(Y)\to X$
of the projections $\pi_t\colon \RipsG_t(Y)\to Y$ for all $t > 0$.
If $s > t$, then both $h_s$ and $h_t$ can be considered as lifts
of $\pi_t\colon \RipsG_t(Y)\to Y$, so they are ls-equivalent.
Notice all $h_t\colon Y\to X$ are bornologous.
Indeed, if $d_Y(x,y) < s$ and $s > t$, then $d_X(h_s(x),h_s(y)) \leq L_s$,
where $L_s$ is the Lipschitz constant of $h_s$.
Since there is $c > 0$ such that $d_X(h_s(z),h_t(z)) < c $ for all $z\in Y$,
$d_X(h_t(x),h_t(y)) \leq L_s+2c$.

It remains to show that $h_s\circ f$ is ls-equivalent to $id_X$
for some $s > 0$. Choose $ s > 0$ such that $d_X(x,y) < 1$
implies $d_Y(f(x),f(y)) < s$ and notice $f\colon \RipsG_1(X)\to \RipsG_s(Y)$
is short.
The map $g=h_s\circ f\colon \RipsG_1(X)\to X$ is Lipschitz
and $f\circ g$ is ls-equivalent to $f\circ \pi_1\colon \RipsG_1(X)\to Y$.
That means $g$ is ls-equivalent to $id_X$ (use the uniqueness of lifts
up to ls-equivalence).
\edokaz

\section{Coarse graphs}

\begin{Def}
Given a graph $G$, by $A(G)$ we mean the graph with the same set of vertices
as $G$ but the set of edges is increased by adding all $[v,w]$
such that $d_G(v,w)=2$.
\end{Def}

In other words, $A(G)$ equals $\RipsG_t(G,d_G)$ for any $2 < t < 3$.

Notice the identity function $G\to A(G)$ is short and bi-Lipschitz for any graph $G$.

\begin{Def}
A {\bf coarse graph}
is a direct sequence $\{V_1\to V_2\to\ldots\}$
of graphs $V_n$ and short maps
$i_{n,m}\colon V_n\to V_m$ for all $n\leq m$ such that

\begin{enumerate}
\item $i_{n,n}=id$ for all $n\ge 1$,
\item $i_{n,k}= i_{m,k}\circ i_{n,m}$ for all $n\leq m\leq k$,
\item for every $n\ge 1$ there is $m > n$ so that $i_{n,m}\colon A(V_n)\to V_{m}$ is short.
\end{enumerate}
\end{Def}

\begin{Def}
A {\bf coarse graph of a metric space $(X,d_X)$}
is a coarse graph $\{V_1\to V_2\to\ldots\}$
together with Lipschitz maps
$p_{n}\colon V_n\to X$
$$
\xymatrix{ V_1 \ar_{p_1}[d] \ar^{i_{1,2}}[r]&    V_2 \ar_{p_2}[dl]   \ar^{i_{2,3}}[r]&V_3   \ar^{p_3}[dll] \ar^{i_{3,4}}[r]&\cdots\\
X &&&
}
$$

for all $n$ such that the following conditions are satisfied:
\begin{itemize}
\item[a.] $p_n$ is ls-equivalent to $p_{n+1}\circ i_{n,n+1}$ for each $n\ge 1$,
\item[b.] For any Lipschitz map $g\colon V\to X$ from a graph $V$
there is $n\ge 1$ and a short map $g^\prime\colon V\to V_n$ such that $p_n\circ g^\prime$
is large scale equivalent to $g$,
\item[c.] If $g,h\colon V\to V_n$ are two short maps from a graph $V$ such that $p_n\circ g$ is large scale equivalent
to $p_n\circ h$, then there is $m > n$ such that
$i_{n,m}\circ g$ is large scale equivalent to $i_{n,m}\circ h$.
\end{itemize}
\end{Def}

Traditionally, a {\bf Cayley graph} of a finitely generated group $G$
is defined to be $\Gamma(G,S)$, where $S$ is a symmetric finite set of generators
of $G$ not containing the neutral element $1_G$.
Its set of vertices is $G$ and edges of $\Gamma(G,S)$
are precisely of the form $[g,g\cdot s]$, $g\in G$ and $s\in S$.
However, one can easily generalize the concept of Cayley graphs
to arbitrary groups $G$ and arbitrary finite subsets  $S$ of $G\setminus \{1_G\}$ (actually, $S$ may be infinite but we do not know any application of such graphs).

It is known that arbitrary countable group $G$ has a proper left invariant metric $d_G$
(see \cite{Sha} or \cite{Smith}) and any two such metrics are coarsely
equivalent. Our next result is a variant of that fact. 

\begin{Prop}\label{CayleyGraphs}
Suppose $d_G$ is a left invariant proper metric on a group $G$.
If $S_n$ is an increasing sequence of finite symmetric subsets of $G$ such that 
$G\setminus \{1_G\}=\bigcup\limits_{n=1}^\infty S_n$, then the sequence
$\{\Gamma(G,S_1)\to \Gamma(G,S_2)\to\ldots\}$ together with projections $\pi_{i}\colon \Gamma(G,S_i)\to G$
forms a coarse graph of $(G,d_G)$.
\end{Prop}

\dokaz Obviously, maps $i_{n,m}\colon \Gamma(G,S_n)\to \Gamma(G,S_m)$
are identities on vertices and are short for $n\leq m$.
For any $n$ there is $m\ge n$ such that $s_n\cdot S_n\subset S_m$
in which case $i_{n,m}\colon A(\Gamma(G,S_n))\to \Gamma(G,S_m)$
is short.
If $g\colon V\to G$ is a short map from a graph $V$ to $G$,
we pick $S_n$ containing the ball $B(1_G,2)$ at $1_G$ of radius $2$
(as $d_G$ is proper, that ball is finite). Now $g$ considered
as a map from $V$ to $\Gamma(G,S_n)$ is short.
Indeed, if $[v,w]$ is an edge in $V$, then $d_G(g(v),g(w))\leq d_V(v,w)=1$,
so $d_G(1_G,g(v)^{-1}\cdot g(w))\leq 1$ and $s=g(v)^{-1}\cdot g(w)\in S_n$.
Therefore $[g(v),g(w)]$ is an edge in $\Gamma(G,S_n)$ proving $g$ is short.
\par Suppose $g,h\colon V\to \Gamma(G,S_n)$ are two short maps
so that $\pi_n\circ g\sim_{ls} \pi_n\circ h$.
There is $M > 0$ such that $d_G(g(v),h(v)) < M$ for all vertices $v$ of $V$.
Choose $m > n$ with the property that $S_m$ contains the ball $B(1_G,M)$.
As above, we can show $[g(v),h(v)]$ is an edge in $\Gamma(G,S_m)$
for any vertex $v$ of $V$. Thus $i_{n,m}\circ g\sim_{ls} i_{n,m}\circ h$.
\edokaz

\begin{Example} The following example generalizes Rips graphs. Suppose $(X,d_X)$ is a metric space and $\mathcal{U}_n$ is a sequence of uniformly bounded covers of $X$ such that 
$\mathcal{U}_{n}$ is a  refinement of $\mathcal{U}_{n+1}$ and Lebesgue numbers
$L(\mathcal{U}_n)$ of $\mathcal{U}_n$  form a  sequence diverging to infinity. The Rips coarse graph of $X$ with respect to the sequence $\mathcal{U}_n$ is $\{\RipsG_{\mathcal{U}_1}(X)\to \RipsG_{\mathcal{U}_2}(X)\to\ldots\}$ together with projections $\pi_{i}\colon \RipsG_{\mathcal{U}_i}(X)\to X$ where
\begin{enumerate}
  \item maps $\pi_i$ are induced by the identity $id_X$;
  \item maps $i_{n,m}\colon \RipsG_{\mathcal{U}_n}(X)\to \RipsG_{\mathcal{U}_m}(X)$ are identities on vertices.
\end{enumerate}

\end{Example}

\begin{Prop}\label{RipsGisCoarse}
If $\mathcal{U}_n$ is a sequence of uniformly bounded covers of $X$ such that $\mathcal{U}_{n}$ is a  refinement of $\mathcal{U}_{n+1}$ and Lebesgue numbers  form a  sequence $L(\mathcal{U}_n)\to\infty$ as $n\to\infty$, then the sequence
$\{\RipsG_{\mathcal{U}_1}(X)\to \RipsG_{\mathcal{U}_2}(X)\to\ldots\}$ together with projections $\pi_{i}\colon \RipsG_{\mathcal{U}_i}(X)\to X$
forms a coarse graph of $(X,d_X)$.
\end{Prop}

\dokaz
We first prove that $\{\RipsG_{\mathcal{U}_1}(X)\to \RipsG_{\mathcal{U}_2}(X)\to\ldots\}$ is a coarse graph. Because $\mathcal{U}_{n}$ is a  refinement of $\mathcal{U}_{n+1}$ all  identity maps $i_{n,m}$ are short. Furthermore because the edges of $A(\RipsG_{\mathcal{U}_n}(X))$ are a subset of edges of $\RipsG_{\mathcal{U}_m}(X)$ provided $L(\mathcal{U}_m)\geq d$ (where $\mathcal{U}_n$ is $d$ bounded) the map $i_{n,m}\colon A(\RipsG_{\mathcal{U}_n}(X)) \to \RipsG_{\mathcal{U}_m}(X)$ is short for every choice of sufficiently large $m > n$.
The maps $\pi_i$ are $a_i-$Lipschitz where $\mathcal{U}_i$ is $a_i-$bounded cover. Also $p_n$ is ls-equivalent to $p_{n+1}\circ i_{n,n+1}$ for each $n\ge 1$.

Let $g\colon V \to X$ be  $a-$Lipschitz map from a graph $V$. Then $g$ induces a short map $g' \colon V \to  \RipsG_{\mathcal{U}_n}(X)$ for every $n$ such that $L(\mathcal{U}_n) \geq a$ and $\pi_{t(n)} \circ g'=g$ holds.

Suppose $g,h\colon V\to \RipsG_{\mathcal{U}_n}(X)$ are two short maps from a graph such that $\pi_{n}\circ g$ is $d-$close to $\pi_{n}\circ h$. Then $i_{n,m}\circ g$ is $1-$close to $i_{n,m}\circ h$ for every $m\ge n$ so that $L(\mathcal{U}_m) \geq d.$
\edokaz

Note that $i_{n,m}\colon A(\RipsG_{\mathcal{U}_n}(X)) \to \RipsG_{\mathcal{U}_m}(X)$ is short if ${\mathcal{U}_n}$ is a star refinement of ${\mathcal{U}_m}$.

\begin{Cor}
If $t(n)\to\infty$ is increasing, then the sequence
$\{\RipsG_{t(1)}(X)\to \RipsG_{t(2)}(X)\to\ldots\}$ together
with projections $\pi_{t(n)}\colon \RipsG_{t(n)}(X)\to X$
forms a coarse graph of $(X,d_X)$.
\end{Cor}

For any two graphs $G_1$ and $G_2$ let $Short(G_1,G_2)_{ls}$
be the set of large scale equivalence classes
of short maps from $G_1$ to $G_2$.

Given a coarse graph $\{V_1\to V_2\to\ldots\}$
and a graph $V$ we consider
the direct limit of $Short(V,V_1)_{ls}\to Short(V,V_2)_{ls}\to\ldots$ and define it as the set of morphisms from $V$ to $\{V_1\to V_2\to\ldots\}$.

The set of morphisms from $\{W_1\to W_2\to\ldots\}$
to $\mathcal{V}=\{V_1\to V_2\to\ldots\}$
is the inverse limit of $\ldots \to Mor(W_n,\mathcal{V})\to\ldots\to Mor(W_1,\mathcal{V})$.

We can restate the above definition of morphisms between coarse graphs as follows. Suppose $\mathcal{V}=\{V_1\stackrel{i_{1,2}}{\to} V_2\stackrel{i_{2,3}}{\to}\ldots\}$ and $\mathcal{W}=\{W_1\stackrel{j_{1,2}}{\to} W_2\stackrel{j_{2,3}}{\to}\ldots\}$ are two coarse graphs. First we consider a {\bf pre-morphism} $F\colon \mathcal{V} \to \mathcal{W}$ that consists of short maps $f_k\colon V_k \to W_{n_F(k)}$ so that for every $k\ge 1$ there is $m\ge n_F(k+1)$ resulting in $j_{n_F(k),m}\circ f_k\sim_{ls} j_{n_F(k+1),m}\circ f_{k+1}\circ i_{k,k+1}$.

$$
\xymatrix{ W_1  \ar^{j_{1,2}}[r]&    W_2   \ar^{j_{2,3}}[r]&  \cdots  \ar^{j}[r]& W_{n_F(1)}\ar^{j}[r] & \cdots \ar^j[r] &W_{n_F(2)}\ar^{j}[r]&\cdots \ar^j[r] & W_{n_F(3)}\cdots\\
V_1 \ar^{f_1}[urrr]   \ar_{i_{1,2}}[r]& V_2 \ar^{f_2}[urrrr]   \ar_{i_{2,3}}[r]& V_3 \ar_{f_3}[urrrrr]  \ar_{i_{3,4}}[r]&\cdots &&&&
}
$$

Two pre-morphisms $F,G\colon \mathcal{V} \to \mathcal{W}$ are considered to be equivalent if for every $k$  there is $m \geq \max\{n_F(k),n_G(k)\}$ so that $j_{n_F(k),m}\circ f_k \sim_{ls} j_{n_G(k),m}\circ g_k$. The sets of equivalence classes of pre-morphisms form the set of morphisms from $\mathcal{V}$ to $\mathcal{W}$.

\begin{Thm}\label{MainThmOnCoarseGraphs}
If $\mathcal{V}=\{V_1\to V_2\to\ldots\}$ is a coarse graph of $(X,d_X)$
and $\mathcal{W}=\{W_1\to W_2\to\ldots\}$ is a coarse graph of $(Y,d_Y)$,
then there is a natural bijection between
bornologous maps from $X$ to $Y$
and morphisms from $\mathcal{V}$ to $\mathcal{W}$.
\end{Thm}
\dokaz We will follow the notation from the diagram below
$$
\xymatrix{ V_1 \ar_{p_1}[d] \ar^{i_{1,2}}[r]&    V_2 \ar_{p_2}[dl]   \ar^{i_{2,3}}[r]&V_3   \ar^{p_3}[dll] \ar^{i_{3,4}}[r]&\cdots&
W_1 \ar_{q_1}[d] \ar^{j_{1,2}}[r]&    W_2 \ar_{q_2}[dl]   \ar^{j_{2,3}}[r]&W_3   \ar^{q_3}[dll] \ar^{j_{3,4}}[r]&\cdots\\
X &&& & Y &&&
}
$$

Suppose $f\colon X \to Y$ is a bornologous map. Notice that every map $f\circ p_i\colon V_i \to Y$ is Lipschitz as $f$ is bornologous and $f\circ p_i$ is defined on a graph.  As $\mathcal{W}$ is a coarse graph of $(Y,d_Y)$,  each $f\circ p_i$ lifts to a short map $f_i \colon V_i \to W_{n(i)}$. Since $\mathcal{V}$ is a coarse graph of $(X,d_X)$, maps $f_i$ form a pre-morphism $F\colon \mathcal{V}\to \mathcal{W}$. Also note that such construction gives us a unique morphism $F$ with the property $f\circ p_* \sim _{ls} q_*\circ F $.

Let $F\colon \mathcal{V}\to \mathcal{W}$ be a pre-morphism. Consider the Rips coarse graph $\{\RipsG_{1}(X)\stackrel{e_{1,2}}{\to} \RipsG_{2}(X)\stackrel{e_{2,3}}{\to}\ldots\}$
of $(X,d_X)$ together with projections $\pi_{n}\colon \RipsG_{n}(X)\to X$. The identity map $1_X \circ \pi_k\colon \RipsG_{k}(X)\to X$ is short and lifts to a short map $g_k \colon \RipsG_{k}(X)\to V_{m(k)}$. Define a map $f \colon X \to Y$ by $f(x):= q_{n(m(1))}\circ f_{m(1)}\circ g_1(x)$. Note that the properties of coarse complexes  and their morphisms imply that all the maps $X \to Y$ defined by  $f(x):= q_{n(m(k))}\circ f_{m(k)}\circ g_k(x)$ are close to each other which means there is exactly one function $f$ factoring over $F$. Also note that maps $q_{n(m(k))}\circ f_{m(k)}\circ g_k\colon \RipsG_{k}(X)\to Y$ are Lipschitz for every choice of $k$ which means that $f$ is bornologous: if $d_X(x,y)<r$ then $d_Y(f(x),f(y))< L+2C$ where $L$ is the Lipschitz constant of the map $q_{n(m(k))}\circ f_{m(k)}\circ g_k$ and $q_{n(m(k))}\circ f_{m(k)}\circ g_k$ is $C$-close to $q_{n(m(1))}\circ f_{m(1)}\circ g_1(x)$. 
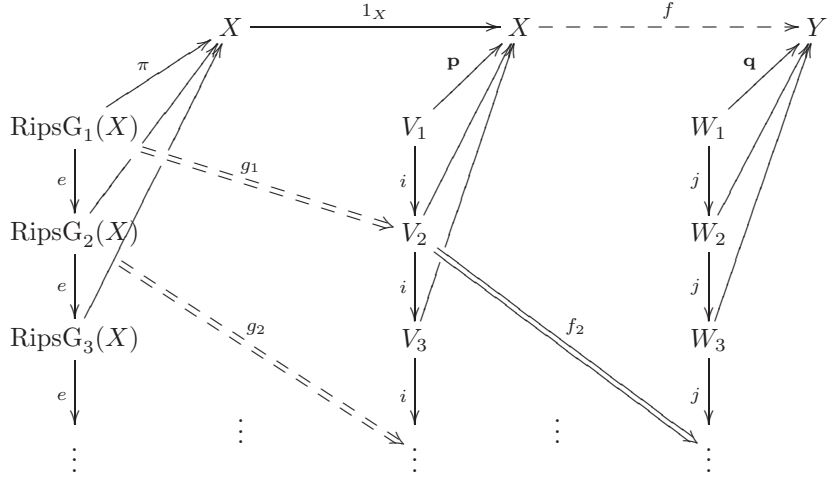
$$
\xymatrix{ & X \ar[rrr]^{1_X} &  &&X \ar@{-->}[rrr]^f &&&Y\\
\RipsG_{1}(X) \ar[d]_e \ar[ur]^\pi \ar@{==>}[rrrd]^{g_1}&&& V_1 \ar[d]_i \ar[ur]^\p&&&W_1 \ar[d]_j \ar[ur]^\q&\\
\RipsG_{2}(X) \ar[d]_e \ar[uur]|(.4)\hole  \ar@{==>}[rrrdd]|(.1)\hole^{g_2} &&&V_2 \ar[d]_i \ar[uur] \ar@{=>}[rrrdd]^{f_2}&&&W_2 \ar[d]_j \ar[uur]&\\
\RipsG_{3}(X) \ar[d]_e \ar[uuur]|(.58)\hole  \ar@{}[rrrd]_{\vdots}&&&V_3\ar@{}[rrrd]_{\vdots}  \ar[d] _i \ar[uuur]|(.27)\hole&&&W_3 \ar[d]_j \ar[uuur]&\\
\vdots &&&\vdots&&&\vdots&
}
$$
The fact that $f$ is a unique map that factors over $f_{m(1)}$ (that is $f\circ p_* \sim _{ls} q_*\circ F $) implies that the rule $F \mapsto f$ is the inverse of the rule $f\mapsto F$ from the previous paragraph. Hence the two rules induce a bijection.

These bijections are natural as $f\circ g \mapsto F\circ G$ which follows from the commutativity of the diagrams.
\edokaz

\begin{Cor} Any two coarse graphs of $X$ are ls-equivalent.
\end{Cor}

\section{Coarse simplicial complexes}

We have seen that graphs are sufficient to describe coarse category
of metric spaces. However, in order to capture more complicated concepts
(asymptotic dimension, coarse connectivity, and Property A) we need to
consider simplicial complexes.

In this section we will consider simplicial complexes $K$ with the set of vertices $X$.
Each such complex induces the graph $G(K)$ obtained by considering only the edges of $K$ (i.e., $G(K)$ is the $1$-skeleton of $K$).

\begin{Prop}
A function $f\colon G(K)\to Y$ is Lipschitz if and only if the family $\{f(\Delta)\}_{\Delta\in K}$
is uniformly bounded in $(Y,d_Y)$.

\end{Prop}
\dokaz If $f$ is $a-$Lipschitz then $\{f(\Delta)\}_{\Delta\in K}$ is uniformly $a-$bounded. Conversely, if $\{f(\Delta)\}_{\Delta\in K}$ is uniformly $a-$bounded then $f$ is $a-$Lipschitz as graphs are $1-$geodesic.
\edokaz

Conversely, each graph $G$ on $X$ induces the minimal flag complex
$F(G)$ containing $G$
(recall $K$ is a {\bf flag complex} if $\Delta$ is a simplex of $K$
whenever $[v,w]$ belongs to $K$ for all $v,w\in\Delta$)
 and each short map
$G_1\to G_2$ induces a simplicial map $F(G_1)\to F(G_2)$.

\begin{Example}
Given a metric space $(X,d)$ and uniformly bounded cover $\mathcal{U}$ of $X$,
{\bf the Rips complex} $\Rips_{\mathcal{U}}(X)$ equals $F(\RipsG_{\mathcal{U}}(X))$.
\end{Example}

We can extend the definition of $A(G)$ from graphs
to complexes as follows:
\begin{Def}
$\Delta\in A(K)$ if and only if there is $v\in K^{(0)}$
such that the set of vertices $\Delta^{(0)}$ of $\Delta$ is contained in the closed
star of $v$ in $K$. Equivalently, there is a vertex $v$ of $K$ such that $[v,w]$
is an edge of $K$ for each vertex $w$ of $\Delta$.

A map $f\colon K \to X$ from a simplicial complex $K$ to a metric space $X$ is \textbf{bornologous} if and only if the family $\{f(\Delta)\}_{\Delta \in K}$ is uniformly bounded.
\end{Def}

Note that any coarse graph $G_1\to G_2\to\ldots$  induces simplicial maps of complexes $F(G_n)\to F(G_m)$.

\begin{Def}
A {\bf coarse simplicial complex}
is a direct sequence $\{K_1\to K_2\to\ldots\}$
of simplicial complexes $K_n$ and simplicial maps
$i_{n,m}\colon K_n\to K_m$ for all $n\leq m$ such that the following conditions are satisfied:
\begin{itemize}
\item[a.] $i_{n,n}=id$ and $i_{n,k}=i_{n,m}\circ i_{m,k}$ for all $n\leq m\leq k$,
\item[b.] for each $n$ there is $m > n$ such that
 $i_{n,m}\colon A(K_n)\to K_m$ is simplicial.
\end{itemize}

\end{Def}

\begin{Example}
For any metric space $(X,d_X)$ and any increasing
sequence $\{t(n)\}$ diverging to infinity,
the sequence of Rips complexes $K_n=\Rips_{t(n)}(X)$
with identity maps $i_{n,m}\colon K_n\to K_m$
for $n \leq m$ forms a coarse simplicial complex.
\end{Example}

\begin{Example} \label{RipsisCoarse} Suppose $(X,d_X)$ is a metric space and $\mathcal{U}_n$ is a sequence of uniformly bounded covers of $X$ such that Lebesgue numbers $L(\mathcal{U}_n)\to\infty$ as $n\to\infty$ and $\mathcal{U}_{n}$ is a  refinement of $\mathcal{U}_{n+1}$. Then the sequence of Rips complexes $K_{\mathcal{U}_{n}}=\Rips_{\mathcal{U}_{n}}(X)$
with identity maps $i_{n,m}\colon K_n\to K_m$
for $n \leq m$ forms a coarse simplicial complex.
Any such coarse complex will be denoted by $\Rips_\ast(X)$
and called {\bf a coarse Rips complex of $X$}.
\end{Example}

\dokaz Similar as first part of \ref{RipsGisCoarse}. Because  $\mathcal{U}_{n}$ is a  refinement of $\mathcal{U}_{n+1}$,  $[v,w]$ is an edge of $\Rips_{\mathcal{U}_{m}}(X)$ for every edge $[v,w]$ of  $\Rips_{\mathcal{U}_{n}}(X)$ with $n\leq m.$ $\Rips_\ast(X)$ being a flag complex of its one skeleton implies that all the maps $i_{n,m}$ are simplicial. Furthermore, edges of $A(\Rips_{\mathcal{U}_{n}}(X))$ form a  subset of edges of  $\Rips_{\mathcal{U}_{m}}(X)$ for every $m$ with $L(\mathcal{U}_m)\geq d$ (where $\mathcal{U}_n$ is $d$ bounded), hence natural map $A(\Rips_{\mathcal{U}_{n}}(X))\to\Rips_{\mathcal{U}_{m}}(X)$ is simplicial.
\edokaz

The following is similar to Roe's concept
of an anti-\v Cech approximation of a metric space $X$:

\begin{Example}
Suppose $(X,d_X)$ is a metric space
and $\mathcal{U}_n$ is a sequence of uniformly bounded covers of $X$
such that $\mathcal{U}_{n+1}$ is a star refinement of $\mathcal{U}_{n}$ for each $n\ge 1$. Then the sequence
$\mathcal{N}(\mathcal{U}_1)\to \mathcal{N}(\mathcal{U}_2)\to\ldots$ of nerves of covers $\mathcal{U}_n$
forms a coarse simplicial complex
if $i_{n,n+1}(U)$ contains the star $st(U,\mathcal{U}_n)$ for each $U\in \mathcal{U}_n$.
Any such coarse complex will be denoted by \v Cech$_\ast(X)$
and called {\bf a coarse \v Cech complex of $X$}.
\end{Example}

\dokaz For  $U\in \mathcal{N}(\mathcal{U}_n)$ define $i_{n,n+1}(U)$ to be any $V\in \mathcal{N}(\mathcal{U}_{n+1})$ that contains the star of $U$. ž
Maps $i_{n,n+1}(U)$ induce maps $i_{n,m}$. To prove that $i_{n,n+1}(U)$ is simplicial consider simplex $[U_1, \ldots, U_k]\in \mathcal{N}(\mathcal{U}_n) $ i.e.   
$U_1\cap \ldots \cap U_k \neq \emptyset$. Then $f(U_1)\cap \ldots \cap f(U_k) \neq \emptyset$ implies $[f(U_1), \ldots, f(U_k)]\in \mathcal{N}(\mathcal{U}_{n+1})$ 
hence map $i_{n,n+1}$ is simplicial $\forall n$.

Suppose each of the sets  $U_1, \ldots, U_k\in \mathcal{U}_{n}$ has nonempty intersection with $U\in \mathcal{U}_{n}$. Then all the sets 
$i_{n,n+1}(U_1), \ldots, i_{n,n+1}(U_k)$ have nonempty intersection (namely contain the set $U$) hence $[i_{n,n+1}(U_1), 
\ldots, i_{n,n+1}(U_k)]\in \mathcal{N}(\mathcal{U}_{n+1})$ which implies that $i_{n,n+1}\colon A(\mathcal{N}(\mathcal{U}_{n}))\to \mathcal{N}(\mathcal{U}_{n+1})$ is simplicial.
\edokaz

\begin{Def} Simplicial maps $f,g\colon K \to L$ between simplicial complexes are \textbf{contiguous} if for every simplex $\Delta$ of $K$, 
$f(\Delta)\cup g(\Delta)$ is contained in some simplex of $L$.
\end{Def}

\begin{Def}

A {\bf coarse complex} of a metric space $(X,d_X)$ is a
coarse complex $\mathcal{K}=\{K_1\to K_2\to\ldots\}$
together with bornologous functions $p_n\colon K_n\to X$, $n\ge 1$,
satisfying conditions
\begin{itemize}
\item[a.] $p_n$ is ls-equivalent to $p_{n+1}\circ i_{n,n+1}$ for each $n\ge 1$,
\item[b.] For each bornologous function $f\colon K\to X$ from a simplicial
complex $K$ to $X$ there is $n\ge 1$ and a simplicial function $g\colon K\to K_n$
such that $f$ is ls-equivalent to $p_n\circ g$,
\item[c.] If $f,g\colon K\to K_n$ are two simplicial functions
so that $p_n\circ f\approx_{ls} p_n\circ g$, then there is $m > n$
such that $i_{n,m}\circ f$ is contiguous to $i_{n,m}\circ g$.
\end{itemize}
\end{Def}

\begin{Example}
Any coarse Rips complex $\Rips_\ast(X)$ together
with identity functions $\Rips_\ast(X)\to X$
forms a coarse complex of $X$,
if the Lebesgue
numbers $L(\mathcal{U}_n)\to\infty$ as $n\to\infty$.
\end{Example}

\dokaz Suppose $\Delta :=[v_1, \ldots, v_k] \in \Rips_{\mathcal{U}_n}(X)$ where $\mathcal {U}_n$ is $d_n$-bounded cover of $X$. 
Then $p_n(\Delta)$ is $d_n$-bounded hence $p_n$ is bornologous. Also note that $p_n=p_{n+1}\circ i_{n,n+1}$ by the definition.

Suppose $f\colon K\to X$ is a bornologous function from a simplicial
complex $K$ to $X$ so that sets $\{f(\Delta)\}_{\Delta \in K}$ are $a-$bounded. Then the naturally induced map $f_n\colon K \to \Rips_{\mathcal{U}_n}(X)$ (naturally induced meaning $p_n\circ f_n=f$) is simplicial for every $n$ with $L(\mathcal{U}_n) \geq a$.

Suppose $f,g\colon K\to K_n$ are two simplicial functions
so that $p_n\circ f$ is $d-$close to $ p_n\circ g$. Then $i_{n,m}\circ f$ is contiguous to $i_{n,m}\circ g$ for every $m\geq n$ with $L(\mathcal{U}_m) \geq d+b$, where  $\mathcal{U}_n$ is $b-$bounded.
\edokaz

\begin{Example}
Any coarse \v Cech complex \v Cech$_\ast(X)$ together
with functions $p_n\colon \mathcal{U}_n\to X$
such that $p_n(U)\in U$ for all $U\in \mathcal{U}_n$
forms a coarse complex of $X$.
\end{Example}
\dokaz
Define $p_n(U)$ to be any point of $U$. If the cover $\mathcal{U}_n$ is $a_n$-bounded then $p_n$ is $2a_n$-bornologous and $p_n$ is $a_{n+1}$ close to $p_{n+1}\circ i_{n,n+1}$.

Suppose $f\colon K \to X$ is a $b$-bornologous map from a simplicial complex $K$ to $X$. Pick $n$ so that $L(\mathcal{U}_n)\geq b$ and define $g\colon K \to \Rips_{\mathcal{U}_n}(X)$ by mapping $v$ to any element of $\mathcal{U}_n$ containing $B(f(v),b)$. Note that $p_n \circ g$ is $a_n$ close to $f$. Furthermore $g$ is simplicial: if $[v_1, \ldots, v_k]\in K$ then $f(v_i)\in g(v_j), \forall i,j$, hence $[g(v_1), \ldots, g(v_k)]\in \Rips_{\mathcal{U}_n}(X)$.

Suppose $g,f\colon K \to \Rips_{\mathcal{U}_n}(X)$ are two simplicial maps, so that $p_n \circ f$ and $p_n \circ g$ are $d-$close and $b-$bornologous.  Choose $m \geq n$ so that $L(\mathcal{U}_m)\geq d+b$ and let $\Delta=[v_1, \ldots, v_k]\in K$. There exists $U \in \mathcal{U}_m$ containing $B(p_m\circ f(v_1),d+b)$ hence it contains $p_m \circ f(v_i), p_m \circ g (v_j), \forall i,j.$ Such set will be contained in all sets $i_{n,m+1}\circ f (v_1)$ and $i_{n,m+1}\circ g (v_1)$ hence $i_{n,m+1}\circ f (\Delta)\cup i_{n,m+1}\circ g (\Delta)$ is contained in a simplex. The parameter $m$ does not depend on the choice of $\Delta$ which implies that $i_{n,m+1}\circ f$ and $i_{n,m+1}\circ g$ are contiguous.
\edokaz

\begin{Prop}
Suppose $L$ is a simplicial complex and
$\mathcal{K}=\{K_1\to K_2\to\ldots\}$ is a coarse simplicial
complex. Consider the direct limit of the direct sequence
of sets of simplicial maps $\{SM(L,K_1)\to SM(L,K_2)\to\ldots\}$
and observe contiguity induces an equivalence relation on that set.
The set of all contiguity classes is the set of ls-morphisms $LS(L,\mathcal{K})$
from
$L$ to $\mathcal{K}$.

If $\mathcal{L}=\{L_1\to L_2\to\ldots\}$ is a coarse complex,
then $\ldots LS(L_n,\mathcal{K})\to LS(L_{n-1},\mathcal{K})\to\ldots\to LS(L_1,\mathcal{K})$ is an inverse sequence and its inverse limit
forms the set of ls-morphisms $LS(\mathcal{L},\mathcal{K})$
from $\mathcal{L}$ to $\mathcal{K}$.
\end{Prop}

\begin{Prop}
If $(X,d_X)$ and $(Y,d_Y)$ are metric spaces
with coarse complexes $\mathcal{K}$ and $\mathcal{L}$ respectively,
then there is a natural bijection between the set of ls-morphisms from $X$ to $Y$
and the set of ls-morphisms from $\mathcal{K}$ to $\mathcal{L}$.
\end{Prop}
\dokaz The proof amounts to a modification of the proof of \ref{MainThmOnCoarseGraphs}.
\edokaz

\begin{Cor}
Any coarse Rips complex of $X$ and any coarse \v Cech complex of $X$ are ls-equivalent.
\end{Cor}

\section{Asymptotic dimension of coarse simplicial complexes}

\begin{Def}
We say that $K\overset{g}\rightarrow M\overset{h}\rightarrow L$ is a {\bf contiguous factorization} of
$K\overset{f}\rightarrow L$ if $f$, $g$, and $h$ are simplicial maps and $f$ is contiguous to $h\circ g$.
\end{Def}

\begin{Def} Given a coarse simplicial complex $\mathcal{K}$
we say its {\bf asymptotic dimension} is at most $n$ (notation:
$\asdim(\mathcal{K})\leq n$) if for each $m$ there is
$k > m$ such that $i_{m,k}$ factors contiguously
through an $n$-dimensional simplicial complex.
\end{Def}

\begin{Cor}
If $\mathcal{K}$ is coarsely equivalent to $\mathcal{L}$, then $\asdim(\mathcal{K})=\asdim(\mathcal{L})$.
\end{Cor}
\dokaz
Let $\mathcal{K}=\{K_1\stackrel{i_{1,2}}{\to}
K_2\stackrel{i_{2,3}}{\to}\ldots\}$ and
$\mathcal{L}=\{L_1\stackrel{j_{1,2}}{\to}
L_2\stackrel{j_{2,3}}{\to}\ldots\}$ be coarsely equivalent.
Let $\varphi\colon \mathcal{K}\to \mathcal{L}$ be an isomorphism and
$\psi\colon \mathcal{L}\to \mathcal{K}$
its inverse. For every $k$ there exist $\alpha(k)$ and $\beta(k)$ such
that the simplicial maps
$\varphi_k\colon K_k\to L_{\alpha(k)}$ and $\psi_k\colon L_k\to
K_{\beta(k)}$ are short.
Because $\psi$ is the inverse of $\varphi$,
$\psi_{\alpha(k)}\varphi_k\sim_{ls} i_{k,\beta(\alpha(k))}$ and
$\varphi_{\beta(k)}\psi_k\sim_{ls} j_{k,\alpha(\beta(k))}$
for large $k$.
Suppose $\asdim(\mathcal{K})\le n$. For every $m$ there exist $k>m$, an $n$-dimensional simplicial complex $M$
and simplicial maps $f\colon K_{\beta(k)}\to M$ and $g\colon M\to K_k$
such that $gf$ and $i_{k,\beta(m)}$ are contiguous.
Let $\widetilde f=f\psi_m$ and $\widetilde g=\varphi_k g$, then
$j_{m,\alpha(k)}$ and $\widetilde g \widetilde f$
are contiguous, hence $\asdim(\mathcal{L})\le n$.
\edokaz

\begin{Thm}
$\asdim(X)=\asdim(\Rips_\ast(X))$ for any metric space $X$.
\end{Thm}
\dokaz Given $t > 0$ consider the cover $\mathcal{U}=\{B(x,t)\}_{x\in X}$
of $X$ by $t$-balls and choose a uniformly bounded cover $\mathcal{V}=\{V_i\}_{i\in J}$
of $X$ together with a function $f\colon X\to J$
such that $B(x,t)\subset V_{f(x)}$ for all $x\in X$
and $f$ factors contiguously through an $n$-dimensional simplicial complex $K$
as $X\overset{g}\to L\overset{h}\to J$, where $L$ is the set of vertices of $K$.

Given $l\in L$ define $W_l$ as the union of all $B(x,t)$ such that $g(x)=l$.
Let us show the multiplicity of $\mathcal{W}=\{W_l\}_{l\in L}$ is at most $n+1$.
Suppose, on the contrary, that
there are mutually different elements $l(0),\ldots,l(n+1)$ of $L$
such that $W_{l(0)}\cap\ldots\cap W_{l(n+1)}\ne\emptyset$.
That means existence of $x\in X$ and elements $x(0),\ldots,x(n+1)$ of $X$
such that $x\in B(x(k),t)$ and $g(x(k))=l(k)$ for all $0\leq k\leq n+1$.
Since $g$ is a simplicial map and $[x(0),\ldots,x(n+1)]$ is a simplex
in the nerve $\mathcal{N}(\mathcal{U})$ of $\mathcal{U}$,
$[l(0),\ldots,l(n+1)]$ is a simplex in $K$ contradicting $K$ being $n$-dimensional.
\par It remains to show $\mathcal{W}$ is uniformly bounded as its Lebesgue number is at least $t$. Given $l\in L$ put $j=h(l)$. If $g(x)=l$ and $h\circ g$ is contiguous to $f$,
$[f(x),j]$ is a simplex in $\mathcal{N}(\mathcal{V})$,
so $V_{f(x)}\cap V_j\ne\emptyset$. That means  $W_l$ is contained in the star
 of $V_{j}$ in $\mathcal{V}$. Therefore
 $\mathcal{W}$ is uniformly bounded.
\edokaz

\section{Connectivity of coarse simplicial complexes}

Here is the basic extension of connectedness to large scale geometry of
metric spaces:
\begin{Def} (\cite{Kap}, Definition 42 on p.19)
A metric space X is {\bf coarsely $k$-connected} if for each $r$ there exists
$R\ge r$ so that the mapping $|\Rips_r(X)|\to |\Rips_R(X)|$ induces a trivial map of 
$\pi_i$ for $0\leq i\leq k$.
\end{Def}

It has a natural generalization to coarse simplicial complexes:

\begin{Def} A coarse simplicial complex
$\mathcal{K}$ is {\bf coarsely $n$-connected}
if for each $m$ there is $k > m$ such that
$i_{k,m}$ induces trivial homomorphisms $\pi_p(i_{k,m})$ of
homotopy groups for all $0\leq p\leq n$.
\end{Def}

Thus $X$ is coarsely $n$-connected
if and only if $\Rips_\ast(X)$
is coarsely $n$-connected.

Corollary 47 of \cite{Kap} says that coarse $k$-connectedness is a quasi-isomorphism invariant.
We can easily generalize it slightly:
\begin{Cor}
Coarse $k$-connectedness is a large scale invariant.
\end{Cor}

Recall that a metric space $(X,d_X)$ is {\bf $t$-chain connected} for some $t > 0$
if for every two points $x,y$ of $X$ there is a $t$-chain joining them
(that means every two consecutive points in the chain are at distance less than $t$).
Alternatively, $\Rips_t(X)$ is connected. Let us show $X$ being coarsely $0$-connected
and $X$ being $t$-chain connected for some $t > 0$ are equivalent concepts.

\begin{Prop}
If $(X,d_X)$ is a metric space,
then the following conditions are equivalent:
\begin{itemize}
\item[a.] $X$ is coarsely $0$-connected,
\item[b.] there is $t > 0$ such that $\Rips_t(X)$ is connected,
\item[c.]  there is $t > 0$ such that $\Rips_s(X)$ is connected for all $s\ge t$,
\item[d.] $d_X$ attains only finite values and there is $t > 0$ such that $H_0(\Rips_t(X))\to H_0(\Rips_s(X))$ is injective for all $s\ge t$.
\end{itemize}
\end{Prop}
\dokaz (a)$\implies$(b). Choose $t > 1$ such that
the image of $\Rips_1(X)$ in $\Rips_t(X)$ is contained in one path-component
of $\Rips_t(X)$. Given two points $x,y\in X$, the corresponding vertices $x$ and $y$
in $\Rips_t(X)$ have to be joinable by a sequence of edges, i.e. $x$ and $y$ are
joinable in $X$ by a $t$-chain.

(b)$\implies$(c). If every two points of $X$ are joinable by a $t$-chain, they
are joinable by $s$-chain for any $s\ge t$ (use the same chain).

(c)$\implies$(d) is obvious.

(d)$\implies$(a). Notice the direct limit of $\tilde H_0(\Rips_n(X))\to \tilde H_0(\Rips_{n+1}(X))\to\ldots$ is trivial and $\tilde H_0(\Rips_t(X))$ maps to it in an injective manner if $n > t$.
Therefore $\tilde H_0(\Rips_t(X))=0$ which means $X$ is $t$-chain connected.
\edokaz

\begin{Def} A metric space $(X,d_X)$
 is {\bf coarsely geodesic}
if it is coarsely equivalent to a geodesic metric space.
\end{Def}

\begin{Prop} Let $(X,d_X)$ be a metric space.
The following conditions are equivalent:
\begin{itemize}
\item [a.] $(X,d_X)$ is coarsely geodesic,
\item[b.] $(X,d_X)$ coarsely $0$-connected and there is  $t > 0$ such that the identity function
$(X,d_X)\to \RipsG_t(X)$ is bornologous,
\item[c.] $(X,d_X)$ coarsely $0$-connected and there is  $t > 0$ such that the identity function
$(X,d_X)\to \RipsG_t(X)$ is a coarse equivalence,
\item[d.] $(X,d_X)$ coarsely $0$-connected and there is  $t > 0$ such that the identity function
$(X,d_X)\to \RipsG_s(X)$ is a coarse equivalence for all $s\ge t$.
\end{itemize}
\end{Prop}
\dokaz If $(Y,d_Y)$ is geodesic, then the identity function
$(Y,d_Y)\to \RipsG_t(Y)$ is large scale Lipschitz for all $t > 0$,
so the identity function
$(X,d_X)\to \RipsG_t(X)$ is bornologous for every coarsely geodesic space $(X,d_X)$
and $t > 0$ sufficiently large.

Suppose the identity function
$(X,d_X)\to \RipsG_t(X)$ is bornologous for some $t > 0$.
Since the identity function
$\RipsG_t(X)\to (X,d_X)$ is $t$-Lipschitz,
both $\RipsG_t(X)$ and $(X,d_X)$ are coarsely equivalent.
Since $\RipsG_t(X)$ is coarsely equivalent to its geometric realization (which is geodesic if $X$ is coarsely $0$-connected and $t$ is sufficiently large),
we are done. \edokaz

\begin{Cor}
Every $t$-geodesic metric space $(X,d_X)$
is coarsely geodesic.
\end{Cor}
\dokaz
Observe $X\to \RipsG_t(X)$ is bornologous. Indeed,
given $x,y\in X$ choose a $t$-chain $x_0,\ldots,x_k$ joining $x$ and $y$
such that $d_X(x,y)=\sum\limits_{i=0}^{k-1} d_X(x_i,x_{i+1})$ and $k$ is minimal
with respect to that property.
Therefore either $d_X(x,y) < 2\cdot t$ or $d_X(x,y) \ge 2\cdot t$, $k\ge 4$,
and $d_X(x_i,x_{i+1})+d_X(x_{i+1},x_{i+2})\ge t$ for each $0\leq i\leq k-2$.
That implies $2\cdot d_X(x,y) \ge 2\cdot (k-1)\cdot t$. Since $k$ is greater than or equal to the distance of $x$ and $y$ in $\RipsG_t(X)$,
 $X\to \RipsG_t(X)$ is bornologous.
\edokaz

\begin{Def} A coarse simplicial complex
$\mathcal{K}$ is {\bf coarsely homology $n$-connected}
if for each $m$ there is $k > m$ such that
$i_{k,m}$ induces trivial homomorphisms $\tilde H_p(i_{k,m})$ of
reduced homology groups for all $0\leq p\leq n$.
\end{Def}

$H_1(X)$ being {\bf uniformly generated} (see \cite{FW}) means there is $R > 0$
such that every loop in $X$ is homologous to the sum of loops,
each of diameter at most $R$.
$\pi_1(X)$ being  {\bf uniformly generated} (see \cite{FW}) means there is $R > 0$ such that
every loop in $X$ extends over a perforated disk $D$ and the image of
each internal hole of $D$ has diameter at most $R$.

Fujiwara and Whyte \cite{FW} proved that every geodesic space $X$
so that $H_1(X)$ is uniformly generated (respectively, $\pi_1(X)$ is uniformly
generated) is quasi-isomorhic isomorphic to a geodesic space $Y$ satisfying
$H_1(Y)=0$ (respectively, $\pi_1(Y)=0$). Their proof involves
adding cones over a family of balls to $X$. Our next result
shows that one can use Rips complexes for the same purpose.

\begin{Prop}
If $X$ is a geodesic metric space,
then the following conditions are equivalent:
\begin{itemize}
\item[a.] $X$ is coarsely homology $1$-connected
(respectively, coarsely $1$-connected),
\item[b.] there is $t > 0$ such that $H_1(\Rips_t(X))=0$
(respectively, $\Rips_t(X)$ is simply connected),
\item[c.]  there is $t > 0$ such that $H_1(\Rips_s(X))=0$
for all $s\ge t$
(respectively, $\Rips_s(X)$ is simply connected for all $s\ge t$),
\item[d.] $H_1(X)$ is uniformly generated
(respectively, $\pi_1(X)$ is uniformly generated).
\end{itemize}
\end{Prop}
\dokaz  If $t < s$, then $H_1(\Rips_t(X))\to H_1(\Rips_s(X))$
is an epimorphism as $X$ is geodesic. That observation takes care of
implications (a)$\implies$(b) and (b)$\implies$(c).

(c)$\implies$(d).
If $H_1(\Rips_t(X))=0$ for some $t > 0$, then any loop can be approximated
by a piecewise-geodesic loop, so it suffices to show that any piecewise-geodesic loop $\gamma$
is homologous to a sum of loops of diameter at most $3t$.
Notice $\gamma$ can be realized in $\Rips_t(X)$, so it is homologous to a sum of
loops representing boundaries of 2-simplices in $\Rips_t(X)$, hence their diameter is at most $3t$. 

(d)$\implies$(c). Consider $t > 0$ such that any element of $H_1(X)$
is homologous to the sum of loops of diameter less than $t$.
Given an element $\gamma$ of $H_1(\Rips_s(X))$ for any $s\ge 3t$ we can realize it as a loop $\gamma\colon S^1\to X$ in $X$ then
extend over an oriented 2-manifold $M$, so that its boundary $\partial M$
is the union $S_0\cup S_1\cup\ldots\cup S_k$, where $S_0=S^1$ and $\gamma(S_i)$ is of diameter
less than $t$ if $i > 0$. We can triangulate $M$ requiring that $\gamma(\Delta)$ has 
diameter less than $t$ for each simplex $\Delta$ of the triangulation.
That triangulation induces a simplicial map from $M$ to $\Rips_s(X)$
that can be extended over cones of each $S_i$, $i > 0$, thus showing that $\gamma$
is homologous to $0$ in $H_1(\Rips_s(X))$.

A similar proof works in the case of the fundamental groups.

\edokaz

\begin{Cor}
If $X$ is a coarsely geodesic metric space,
then $X$ is coarsely equivalent to a simply connected geodesic space
if and only if $X$ is coarsely $1$-connected.
\end{Cor}
\dokaz
By definition $X$ is coarsely equivalent to a geodesic metric space $Y$.

($\Rightarrow$) By the above proposition (see implication d. $\Rightarrow$ a.), a geodesic
1-connected metric space is also
coarsely 1-connected.

($\Leftarrow$) By the above proposition (see implication a. $\Rightarrow$ b.) $Y$ is
coarsely equivalent to the 1-connected
geodesic space $\Rips_t(Y)$ for some $t$.
\edokaz

\section{Coarse trees}

The purpose of this section is to provide a simple proof of a result
of Fujiwara and Whyte \cite{FW} (see \ref{CoarseTreesChar}).

\begin{Thm}
If $X$ is a coarsely geodesic metric space, then the following conditions are equivalent:
\begin{itemize}
\item[a.] $X$ is coarsely equivalent to a simplicial tree,
\item[b.] $\asdim(X)\leq 1$ and $X$ is coarsely
homology $1$-connected,
\item[c.] $\asdim(X)\leq 1$ and $X$ is coarsely
$1$-connected.
\end{itemize}
\end{Thm}
\dokaz  Assume $(X,d_X)$ is geodesic.
\par
(b)$\implies$(c). There is $s > t$ such that $H_1(\Rips_t(X))=0$ and $\Rips_t(X)\to \Rips_s(X)$ contiguously factors through a $1$-complex $K$. 
Since the image of $\pi_1(\Rips_t(X))$ in $\pi_1(K)$ is both perfect and free, it must be trivial.
\par (c)$\implies$(a). Pick a contiguous factorization
$\Rips_t(X)\overset{f}\rightarrow K\overset{g}\rightarrow \Rips_s(X)$ such that $K$ is a simplicial tree and both projections $\pi_t\colon \RipsG_t(X)\to X$
and $\pi_s\colon \RipsG_s(X)\to X$ are coarse equivalences. Such a factorization exists as we can replace $K$ by its universal cover.
We may assume $K=f(\RipsG_t(X))$ as $\RipsG_t(X)$ is connected and we will use
the same notation for $f$ considered as a function defined on $X$ and $g\colon K\to X$.
Therefore $f$ and $g$ are bornologous and $g\circ f\sim_{ls} id_X$.
 As $f$ is bornologous,
there is $M > 0$ such that $d_K(f\circ g\circ f(x),f(x)) < M$ for all $x\in X$.
Given a vertex $v$ of $K$ pick $x\in X$ so that $v=f(x)$. Now $d_K(f\circ g(v),v) < M$
proving $f\circ g\sim_{ls} id_K$.
\edokaz

\begin{Cor}[Fujiwara and Whyte \cite{FW}]\label{CoarseTreesChar}
Suppose $X$ is a geodesic metric space. $X$ is quasi-isometric to a simplicial tree
if $H_1(X)$ is uniformly generated and $X$ is of asymptotic dimension $1$.
\end{Cor}

\ref{CoarseTreesChar} was used in \cite{FW} to show that finitely
presented groups of asymptotic dimension $1$ are virtually free
(see also \cite{Gen} and \cite{JS}).

\section{Property A}

Property A of G.Yu is usually defined for metric spaces of bounded geometry
(that means the number of points in each $r$-ball $B(x,r)$ does not exceed $n(r) <\infty$
for each $r > 0$) as the condition that for each $R,\epsilon > 0$
there is $S > 0$ and finite subsets $A_x$ of $X\times N$, $x\in X$,
so that $A_x\subset B(x,S)\times N$ for each $x\in X$
and $\frac{|A_x\Delta A_y|}{|A_x\cap A_y|} < \epsilon$ if $d(x,y) < R$.
Here $A\Delta B:= (A\setminus B)\cup (B\setminus A)$ is the {\bf symmetric difference}
of sets $A$ and $B$.

For arbitrary metric spaces $X$ one can use
Condition 2 of Theorem 1.2.4 of \cite{Willett}:
\begin{Def}\label{PropertyADef}
 $X$
has {\bf Property A} if and only if for each $R,\epsilon > 0$
there is a function $\xi\colon X\to l^1(X)$ and $S > R$
such that $||\xi_x||_1=1$ for each $x\in X$,
$||\xi_x-\xi_y|| < \epsilon $ if $d(x,y) \leq R$, and $\xi_x$
is supported in $B(x,S)$ for each $x\in X$.
\end{Def}
 
Notice one can always assume $\xi_x$ has non-negative values
(replace $\xi_x$ by its absolute value).
The conditions in \ref{PropertyADef} are weaker than the original
definition of Yu stated in the beginning of this section but both
are equivalent for spaces of bounded geometry (see \cite{Willett}).

In this section we redefine Property A of Yu in terms of coarse
simplicial complexes. Given a simplicial complex $K$
by $\vert K\vert_m$ we mean the geometric realization of $K$
equipped with the metric resulting from considering $\vert K\vert $
as a subset of $l^1(K)$. It is obvious every simplicial map
$f\colon K\to L$ induces a short map $f\colon \vert K\vert_m\to \vert L\vert_m$.

Given two functions $f,g\colon \vert K\vert \to \vert L\vert$
we say they are contiguous if for every simplex $\Delta$ of $K$
there is a simplex $\Delta_1$ of $L$
such that $f(\vert\Delta\vert)\cup g(\vert\Delta\vert)\subset \vert \Delta_1\vert$.
That generalizes the concept of contiguity between simplicial maps of
simplicial complexes.

\begin{Def}\label{PropertyAForCoarseSimplicialComplexes}
A coarse simplicial complex $\mathcal{K}=\{K_1\to K_2\to\ldots\}$
has {\bf Property A} if for each $k\ge 1$ and each $\epsilon > 0$
there is $n > k$ and a function $f\colon \vert K_k\vert_m \to \vert K_n\vert_m$
such that $f$ is contigous to $i_{k,n}\colon \vert K_k\vert\to \vert K_n\vert$
and the diameter of $f(\vert\Delta\vert)$ is at most $\epsilon$ for each
simplex $\Delta$ of $K_k$.
\end{Def}

\begin{Thm}\label{MainThmOnPropertyA}
A metric space $X$ has Property A if and only if its Rips complex has Property A.
\end{Thm}
\dokaz Suppose $X$ has Property A as in \ref{PropertyADef} and $R,\epsilon > 0$.
There is a function $\xi\colon X\to l^1(X)$ and $S > R$
such that $||\xi_x||_1=1$ for each $x\in X$, $\xi_x$ has non-negative values,
$||\xi_x-\xi_y|| < \epsilon $ if $d(x,y) \leq R$, and $\xi_x$
is supported in $B(x,S)$ for each $x\in X$.
By adjusting the values of $\xi_x$ we may assume its support is finite
for each $x\in X$ (pick a finite subset $C(x)$ so that $\sum\limits_{y\in C(x)}\xi_x(y) > 1 -\epsilon/2$
and shift the sum of remaining values to a point in $C(x)$). Therefore $\xi$ may be viewed as
a function from vertices of $\Rips_R(X)$ to $\vert\Rips_{2S}(X)\vert_m$
and can be extended over $\vert\Rips_R(X)\vert_m$ (use the convex structure of $l^1(X)$
and extend the function from vertices over simplices via linear combinations) so that
the resulting $\xi\colon \vert\Rips_R(X)\vert_m\to \vert\Rips_{2S+R}(X)\vert_m$
is contiguous to the inclusion-induced $\vert\Rips_R(X)\vert_m\to \vert\Rips_{2S+R}(X)\vert_m$ and $\xi(\vert\Delta\vert)$ is of diameter at most $\epsilon$ for any simplex $\Delta$ of $\Rips_R(X)$.
\par Conversely, if $\xi\colon \vert\Rips_R(X)\vert_m\to \vert\Rips_{S}(X)\vert_m$
is contiguous to the inclusion-induced $\vert\Rips_R(X)\vert_m\to \vert\Rips_{S}(X)\vert_m$ and $\xi(\vert\Delta\vert)$ is of diameter at most $\epsilon$ for any simplex $\Delta$ of $\Rips_R(X)$,
then $\xi$ restricted to vertices of $\Rips_R(X)$ gives a function $\mu\colon X\to l^1(X)$
by $\xi(x)=\sum\limits_{y\in X}\mu_x(y)\cdot y$. Since for each $x\in X$ there is
a simplex $\Delta_x$ of $\Rips_S(X)$ containing both $x$ and $\xi(x)$,
the carrier of $\mu_x$ is contained in $B(x,S)$. Also, if $d(x,y) < R$, $[x,y]$ forms a simplex
in $\Rips_R(X)$ and its image is of diameter at most $\epsilon$
resulting in $||\mu_x-\mu_y|| \leq\epsilon$.
\edokaz

Since having Property A is an invariant of large scale equivalence of coarse simplicial
complexes, Theorem \ref{MainThmOnPropertyA} really says $X$ has Property A
if and only if any of its coarse simplicial complexes has Property A.

\end{document}